\documentclass[final,1p,times,authoryear,sort]{elsarticle}
\usepackage{amsmath,amsfonts,amssymb,amsthm}
\usepackage[all]{xy}
\usepackage{graphicx}
\usepackage{aliascnt,hyperref}

\journal{arXiv}

\newtheorem{teo}{Theorem}

\def\partet#1#2#3#4{\newaliascnt{#1}{#2}\newtheorem{#1}[#1]{#3}\aliascntresetthe{#1}\providecommand*{#4}{#3}}
\def\parted#1#2#3#4{\newaliascnt{#1}{#2}\newdefinition{#1}[#1]{#3}\aliascntresetthe{#1}\providecommand*{#4}{#3}}
\partet{lema}{teo}{Lemma}{\lemaautorefname}
\partet{prop}{teo}{Proposition}{\propautorefname}
\partet{cor}{teo}{Corollary}{\corautorefname}
\parted{defs}{teo}{Definition}{\defsautorefname}
\parted{ejemplo}{teo}{Example}{\ejemploautorefname}
\parted{obs}{teo}{Remark}{\obsautorefname}
\parted{parr}{teo}{}{\parrautorefname}
\parted{notc}{teo}{Notation}{\notcautorefname}
\parted{ass}{teo}{Assumption}{\assautorefname}
\newproof{dem}{Proof}
\newtheorem*{nteo}{Theorem}

\newtheorem*{nprop}{Proposition}
\newtheorem*{nconj}{Conjecture}
\theoremstyle{definition}
\newtheorem*{ndefs}{Definition}
\newtheorem*{nnotc}{Notation}
\newtheorem*{nobs}{Remark}
\newtheorem*{nejemplo}{Example}

\def\Z{\mathbb{Z}}
\def\C{\mathbb{C}}
\def\N{\mathbb{N}}

\def\PP{\mathbb{P}}
\def\R{\mathbb{R}}

\def\then{\Longrightarrow}
\def\iff{\Longleftrightarrow}

\def\fl#1{\stackrel{#1}{\longrightarrow}}

\def\rk{\text{rk}}
\def\im{\mathrm{Im}}

\begin{document}

\begin{frontmatter}

\title{Solving a sparse system using linear algebra.}

\author[dm]{C\'esar Massri\corref{correspondencia}\fnref{financiado}}
\ead{cmassri@dm.uba.ar}
\address[dm]{Department of Mathematics, FCEN, University of Buenos Aires, Argentina}
\cortext[correspondencia]{Address for correspondence: Department of Mathematics, FCEN,
University of Buenos Aires, Argentina. Postal address: 1428. Phone number: 54-11-4576-3335.}
\fntext[financiado]{The author was fully supported by CONICET, Argentina}

\begin{abstract}
We give a new theoretical tool to solve sparse systems with finitely many solutions.
It is based on toric varieties and basic linear algebra; eigenvalues, eigenvectors and coefficient matrices.
We adapt Eigenvalue theorem and Eigenvector theorem to work with a canonical rectangular matrix (the first
Koszul map) and prove that these new theorems serve to solve overdetermined
sparse systems and to count the expected number of solutions.
\end{abstract}

\end{frontmatter}

\section*{Introduction.}
\subsection*{Overview of the problem.}
In this article we generalize two methods to solve systems of polynomial equations using a coefficient matrix.
One method is based on
the eigenvalue theorem, first noticed in \citet{MR619687}. Another, on the eigenvector theorem, first described in \citet{MR1022943}.
Let us start describing them.

For simplicity, consider a generic system of $n$ polynomial equations with finitely many solutions in $\C^{n}$, all with multiplicity one,
$$\left\{\begin{array}{rcl}
f_1(x_1,\ldots,x_n)&=&0\\
&\vdots&\\
f_n(x_1,\ldots,x_n)&=&0\\
\end{array}\right.$$
where $f_1,\ldots,f_n$ are polynomials in $\C[x_1,\ldots,x_n]$. The quotient ring,
$$\mathcal{R}=\C[x_1,\ldots,x_n]/\langle f_1,\ldots,f_n\rangle,$$
is a finite-dimensional vector space and
its dimension is the number of solutions (we are assuming that
all the solutions have multiplicity one).

Every polynomial $f\in\C[x_1,\ldots,x_n]$, determines a linear map $M_f:\mathcal{R}\rightarrow\mathcal{R}$,
$$M_f(\overline{g})=\overline{fg},\quad g\in\C[x_1,\ldots,x_n],$$
where $\overline{g}$ denotes the class of the polynomial $g$ in the quotient ring $\mathcal{R}$.
The matrix of $M_f$ is called the \emph{multiplication matrix} associated to the polynomial $f$.
\begin{nteo}[Eigenvalue Theorem]
The eigenvalues of $M_f$ are $\{f(\xi_1),\ldots,f(\xi_r)\}$, where $\{\xi_1,\ldots,\xi_r\}$
are the solutions of the system of polynomial equations.
See \citet[Theorem 2.1.4]{MR2161984} for a proof.
\end{nteo}

\begin{nteo}[Eigenvector Theorem]
Let $f=\alpha_1x_1+\ldots+\alpha_nx_n$ be a generic linear form and let $M_f$ be its multiplication matrix.
Assume that $B=\{1,x_1,\ldots,x_n,\ldots\}$ is a finite basis of $\mathcal{R}$ formed by monomials.
Then the left eigenvectors of $M_f$ determine all the solutions of the
system of polynomial equations. Specifically, if $v=(v_0,\ldots,v_n,\ldots)$
is a left eigenvector of $M_f$ such that $v_0=1$, then $(v_1,\ldots,v_n)$ is a solution
of the system of polynomial equations.
See \citet[\S 2.1.3]{MR2161984} for a proof.
\end{nteo}

Now, let us describe the
construction of the coefficient matrix (also in the case of polynomial equations).

Let $d=d_1+\ldots+d_n-n+1$, where $d_i=\deg(f_i)$, $1\leq i\leq n$. Let $S_d$ be the space of
polynomials of degree $\leq d$. Consider the following sets of monomials,
$$\begin{array}{rcl}
B_n&=&\{x_1^{m_1}\ldots x_n^{m_n}\in S_d\,\colon\,d_n\leq m_n\}\\
B_{n-1}&=&\{x_1^{m_1}\ldots x_n^{m_n}\in S_d\setminus B_n\,\colon\,d_{n-1}\leq m_{n-1}\}\\
&\vdots&\\
B_{1}&=&\{x_1^{m_1}\ldots x_n^{m_n}\in S_d\setminus B_2\,\colon\, d_1\leq m_1\}\\
B_0&=&\{x_1^{m_1}\ldots x_n^{m_n}\in S_d\setminus B_1\}.
\end{array}$$
Using these sets, we can consider the following linear map,
$$\Psi:\langle B_0\rangle\times \ldots\times \langle B_n\rangle\rightarrow S_d,\quad \Psi(g_0,\ldots,g_n)=
f_0\cdot g_0 + \sum_{i=1}^n f_i\cdot g_i,$$
where the polynomial $f_0$ is a generic linear form and $\langle B_i\rangle$ is
the vector space generated by $B_i$, $0\leq i\leq n$.
The \emph{coefficient matrix} $M$ is the matrix of $\Psi$ in the monomial bases $B_0,\ldots,B_n$.
It is a square matrix and can be divided into four blocks,
$$M=\begin{pmatrix}M_{11}&M_{12}\\M_{21}&M_{22}\end{pmatrix}.$$

The relation between the coefficient matrix and the multiplication matrix is the following,
\begin{nteo}
For generic systems $f_1,\ldots,f_n$ in $n$ variables, the multiplication matrix associated to $f_0$
in $\mathcal{R}$ is the Schur complement of $M_{22}$ in the coefficient matrix $M$,
$$M_{f_0}=M_{11}-M_{12}M_{22}^{-1} M_{21}.$$
See \citet{Emiris:1994:MBP:190347.190374} and \citet{MR1762401} for a proof.
\end{nteo}

There are several technical difficulties in order to generalize the previous constructions.
For example, the choices of the sets $B_0,\ldots,B_n$ and
the fact that we need a generic system of $n$ polynomial equations in $n$ variables.
The sets $B_0,\ldots,B_n$ are given to assure that
$M_{22}$ is a non-degenerate matrix and that $M$ is a square matrix.
Another technical difficulty is that the system must have simple roots and that $f_0$ must be linear.
All these difficulties may be solved to give generalizations of the constructions,
not only to polynomial equations, but also to sparse systems.
See in the next subsection for the existing work.

In this article, we propose a simpler approach to deal not only with polynomial
equations, but also with sparse systems in general.
We make a canonical choice for the map $\Psi$ (the first Koszul map)
and we make no assumption on $f_0$ nor on the multiplicities of the solutions.
We construct a matrix $M_{11}+M_{12}F$, where $F$ satisfies the linear equation $M_{22}F=-M_{21}$
and such that every solution $\xi$ of the sparse system determines an eigenvalue $f_0(\xi)$
and a left eigenvector of $M_{11}+M_{12}F$. The matrix $M_{11}+M_{12}F$ can be obtained
by elementary column operations on $M$.

Our construction can be used to solve overdetermined sparse systems and also,
to count the expected number of solutions.
The main problem of our matrix $M$ is its size.

\subsection*{Existing work.}
Several classes of scientific and engineering problems are expected to reduce
to algebraic systems with sparse structure. Sparse systems are typical for such a situation.
For example, problems in vision \citep{emiris:inria-00073580},
edge detection,
robot kinematics (kinematics of molecueles/mechanisms), calibration of Gough/Stewart platforms
\citep{933121,Mourrain:1993:LPP:164081.164120}, structural biology and computational chemistry
\citep{MR1703585}.

Given a sparse system, we could ask if there exist solutions. Just as in the affine case,
where the classical Hilbert Nullstellensatz is available, we can apply the Sparse
Nullstellensatz to obtain an answer \citep[Theorem 2.13]{MR1659402}
\begin{nteo}[Corollary of the Sparse Nullstellensatz]
If the ideal generated by $f_1,\ldots,f_k$ contains the unity
$1\in\langle f_1,\ldots,f_k\rangle$, then
the sparse system has no solution in $(\C\setminus 0)^{n}$.
\end{nteo}
The most common way to check the hypothesis of this theorem is using elimination theory \citep{MR1142904}.
The central object in elimination theory is the resultant, which characterizes the solvability of a sparse system
with prescribed support.
The resultant is a polynomial in the coefficient of the sparse system, $\{f_1,\ldots,f_n\}$.
It provides a necessary (and generically sufficient) condition for the existence of solutions.
If the system has a solution, the resultant $R_{f_1,\ldots,f_n}$ is non-zero.
The most famous example of resultant is the determinant of a system of linear equations.

The first mathematicians who worked in elimination theory were Gauss, B\'ezout and Euler in the eighteenth century.
The study of resultants, in the second half of the nineteenth century, started with Sylvester, Cayley, Macaulay and Dixon.
In the last decade of the twentieth century, the theory was reborn with the pioneering work of Jouanolou in 1991 \citep{MR1142904}.
Today, the resultant may be considered, not only in affine or projective space, but also in the toric case.
The foundations were laid in the work of Gelfand,
Kapranov and Zelevinsky \citep{MR1264417}.
Subsequence papers extended the theory into several different directions, see \citet{MR1174606,MR1245200}.

In order to compute the resultant, several algorithm are given.
In \citet{MR1251972} the authors proposed a
formula for the resultant of a system of $n+1$ Laurent polynomials in $n$ variables. They constructed a matrix
whose determinant is a non-zero multiple of the resultant.
This construction is closely related to that of Macaulay's, who called these matrices, \emph{coefficient matrices}, \citep{MR1577000}.
In general, the construction of the coefficient matrices needs a clever choice of monomials.
In \citet{Emiris:1994:MBP:190347.190374}
the authors used a coefficient matrix to obtain a monomial basis for the coordinate ring generated by the given polynomials.

A number of methods exist for constructing matrices whose determinant is the resultant
or, more generally, a non-trivial multiple of it.
These matrices represent the most efficient way for computing
the resultant and for solving sparse systems by means of the resultant method.
For the classical resultant method see \citet{MR619687,MR1056626}.
For the sparse resultant see \citet{Emiris94sparseelimination}, where
an efficient and general algorithm is given. The author studied the complexity
of the algorithm and also the numerical issues.

There are several articles that used
coefficient matrices and/or multiplication matrices to compute the solutions of a system of polynomial equations.
For example, in \citet{MR1022943,MR1799632} and \citet{MR2391492}
the authors gave an algorithm to compute the solutions of
a system of polynomial equations with the same number of variables and equations.
In \citet[\S 6.2]{MR2322065}, the authors showed a generalization of the method to solve an overdetermined system
of polynomial equations
and in \citet{emiris:inria-00073580,MR1398322,MR1832413,MR1709416} and \citet{Emiris:1994:MBP:190347.190374}
the authors gave another generalization, but
to solve a sparse system with the same number of variables and equations.

In \citet{MR1374227}, the authors gave an algorithm using a coefficient matrix
that can treat an overdetermined sparse system.
The authors wrote
\emph{``An important aspect of the algorithm is that it readily extends to
systems of more than $n+1$ polynomials in $n$ variables''}.
They proposed a method to construct the coefficient matrix minimizing its
size. This method was implemented in \citet{emiris:inria-00073580}.

As a final remark, let us mention that there exists another theory
to solve a system of equations using a topological point of view.
It is a called \emph{homotopy method}.
Essentially, first define a trivial system of equations to which all
solutions are easily known. Then, deform the trivial system into the original system.
As the system is deformed the solutions are deformed also, thereby creating paths of solutions.
These paths start from each of the trivial solutions and connect to the solutions
of the original system. By following these paths from the trivial system,
all the solutions of the original system can be determined, see \citet{MR914807}.

\subsection*{Main result.}
We propose a general framework to solve a sparse system using a rectangular coefficient matrix.
Known methods require the construction of a square matrix adapted to each specific system, see for example \citet[3.1]{MR1799632} and \citet[\S 3.2.3]{MR1662032}. One advantage of this new method
is that the construction of the rectangular matrix is \emph{canonical} and does not require a clever
choice of the monomials for its construction.
Our contribution to the theory is the exposure of the properties of the rectangular coefficient matrix $M$
associated to the first Koszul map of a sparse system.

Given that our coefficient matrix is rectangular, it is not possible to
use the previous theorems where a square matrix is required (see \cite{MR1398322} for the sparse case).
Hence, we adapted them to our requirements. This means that we generalized known
theorems to the case of a rectangular coefficient matrix.

\bigskip

Let us list the main results of this article (for definitions and notations see below).
Let $f_0,\ldots,f_k$ be Laurent polynomials
with Newton polytopes $\mathcal{A}_0,\ldots,\mathcal{A}_k$ respectively.
Let $\mathcal{B}_i=\mathcal{A}_0+\ldots+\widehat{\mathcal{A}_i}+\ldots+\mathcal{A}_k$, $0\leq i\leq k$
and let $\mathcal{E}=\mathcal{A}_0+\ldots+\mathcal{A}_k$.

The coefficient matrix $M$ associated to the sparse system $\{f_1,\ldots,f_k\}$ and $f_0$ is
the matrix of $\Psi$ in the monomial bases $\mathcal{B}_i\cap\Z^n$, $0\leq i\leq k$ and $\mathcal{E}\cap\Z^n$,
$$\Psi:S_{\mathcal{B}_0}\times\ldots\times S_{\mathcal{B}_k}\fl{} S_{\mathcal{E}},\quad
\Psi(g_0,\ldots,g_k)=f_0\cdot g_0+\sum_{i=1}^k f_i\cdot g_i.$$
Matrix $M$ is rectangular and can be divided into four blocks,
$$M=\begin{pmatrix}
M_{11}&M_{12}\\
M_{21}&M_{22}
\end{pmatrix},\quad M_{11}\in\C^{p\times p},\,p=\dim(S_{\mathcal{B}_0}).$$

\bigskip

\noindent {\bf Main Hypotheses:}
Assume that $0\in \mathcal{A}_0$ and $f_0$ is a non-constant Laurent polynomial,
that the lattice polytope $\mathcal{E}$ is full dimensional
and finally, that $\langle f_0,f_1,\ldots,f_k\rangle=S$.
\hspace{\fill}\qed

\bigskip

Using the matrix $M$, we can test the last assumption adapting a theorem due to Macaulay, see \citet{MR1577000}
or \citet[Theorem 3.7]{MR1662032}.
\begin{nprop}[\autoref{spar-null}]
Assume that $\mathcal{E}$ is full dimensional.
Then, $M$ has full rank if and only if $\langle f_0,\ldots,f_k\rangle=S$.
\hspace{\fill}\qed
\end{nprop}

Using the previous new proposition and as a benefit of our approach,
we obtained a proof of a conjecture due to J. Canny and I. Emiris, \citep{MR1768142}.
\begin{nconj}[8.3, Sparse Effective Nullstellensatz over $\C$]
Suppose $f_0,\ldots,f_k$ are arbitrary Laurent polynomials in $S=\C[x_1^{\pm1},\ldots,x_n^{\pm1}]$
with Newton polytopes $\mathcal{A}_i$, $0\leq i\leq k$ such that
the generated ideal is $S$, $\langle f_0,\ldots,f_k\rangle =S$.
Then there exist Laurent polynomials $g_0,\ldots,g_k\in S$, with Newton polytopes $\mathcal{B}_i$, $0\leq i\leq k$,
such that
$$1=\sum_{i=0}^kf_i\cdot g_i,\quad
\mathcal{B}_i\subseteq\mathcal{A}_0+\ldots+\widehat{\mathcal{A}_{i}}+\ldots+\mathcal{A}_k.$$
\end{nconj}
\begin{dem}
Given that $\langle f_0,\ldots,f_k\rangle =S$, $\Psi$ is surjective, hence $1\in S_{\mathcal{E}}$.
\hspace{\fill}\qed
\end{dem}

Another new result that we proved is
a formula to count the number of expected solutions of a sparse system using $M$
and also, our main theorem; an adaptation of the Eigenvalue/Eigenvector Theorem to the case
of a rectangular coefficient matrix.
\begin{nteo}[\autoref{spar-teo1}(a)]
The sparse system $\{f_1,\ldots,f_k\}$ has a finite number of expected solutions equal to
$$\rk(M)-
\rk
\begin{pmatrix}
M_{12}\\
M_{22}
\end{pmatrix}\geq 0.$$
\end{nteo}

\begin{nteo}[\autoref{spar-F},\autoref{spar-mult}]
Let $F$ be a solution of the linear equation $M_{21}+M_{22}F=0$.
Then, every solution, $\xi$, of the sparse system determines a
left eigenvector of $M_{11}+M_{12}F$ with eigenvalue $f_0(\xi)$.
The multiplicity of $f_0(\xi)$ is greater than or equal to the multiplicity of $\xi$.
\hspace{\fill}\qed
\end{nteo}

\subsection*{Summary.}
This paper is organized as follows.
In {\bf \autoref{sec-prel}} we present some preliminaries about toric varieties and lattice polytopes.
In \autoref{prel-mink} we construct an irreducible projective toric variety $X$ associated to a full dimensional
lattice polytope $\mathcal{E}$
and relate $\N$-Minkowski summands of $\mathcal{E}$ with invertible sheaves on $X$ generated by their global sections.
In {\bf \autoref{sec-spar}} we construct a stably twisted Koszul complex and we apply it in
two different ways.
Firstly, in \autoref{spar-F}, we use it to prove
that every solution of the sparse system determines a left eigenvector/eigenvalue of a matrix built from this complex.
Secondly, in \autoref{spar-teo1}(a), we use it to
count the number of expected solutions of the sparse system (counted with
multiplicities).
In {\bf \autoref{sec-app}} we give an application.

\section{Preliminaries.}\label{sec-prel}
A \emph{sparse system} is a collection of Laurent polynomials, $\{f_1,\ldots,f_k\}$,
$$f_i=\sum_{v\in\mathcal{Q}_i}c_{i,v}x_1^{v_1}\ldots x_n^{v_n},\quad 1\leq i\leq k,$$
where $\mathcal{Q}_i$ are fixed finite subsets of $\Z^n$.
The set $\mathcal{Q}_i$ is called the \emph{support} of $f_i$.
The convex hull $\mathcal{A}_i$ of $\mathcal{Q}_i$,
$$\mathcal{A}_i=\text{conv}(\mathcal{Q}_i)\subseteq\R^n,$$
is called the \emph{Newton polytope} of $f_i$, denoted $N(f_i)$, $1\leq i\leq k$.

\begin{ndefs}
A \emph{lattice polytope} $\mathcal{A}\subseteq\R^n$ is the convex hull of a finite set $\mathcal{Q}\subseteq\Z^n$,
$\mathcal{A}=\text{conv}(\mathcal{Q})$.

The \emph{dimension} of a lattice polytope $\mathcal{A}\subseteq\R^n$, is the dimension of the smallest
affine subspace of $\R^n$ containing $\mathcal{A}$. We say that $\mathcal{A}$ is a \emph{full dimensional}
lattice polytope when the dimension of $\mathcal{A}\subseteq\R^n$ is $n$.
\end{ndefs}

\begin{nnotc}
Let $S=\C[x_1^{\pm1},\ldots,x_n^{\pm1}]$ be the algebra of Laurent polynomials.
Given a lattice polytope $\mathcal{A}$, let $S_{\mathcal{A}}$ be the vector space of polynomials
with Newton polytopes in $\mathcal{A}$,
$$S_{\mathcal{A}}=\{g\in S\,\colon\,N(g)\subseteq\mathcal{A}\}.$$
The dimension of $S_{\mathcal{A}}$ is equal to the cardinal of $\mathcal{A}\cap\Z^n$,
$$\dim(S_{\mathcal{A}})=\#(\mathcal{A}\cap\Z^n).$$
The finite set $\mathcal{A}\cap\Z^n$ determines a monomial basis for $S_{\mathcal{A}}$.
\end{nnotc}

\begin{ndefs}
Given lattice polytopes $\mathcal{B}$ and $\mathcal{E}$ in $\R^n$,
we say that $\mathcal{B}$ is an \emph{$\N$-Minkowski summand} of $\mathcal{E}$
if
$$\mathcal{B}+\mathcal{B}'=k\mathcal{E}$$
for some positive integer $k$ and lattice polytope $\mathcal{B}'\subseteq\R^n$.

For example, $2\mathcal{E}$ is an $\N$-Minkowski summand of $\mathcal{E}$.
\end{ndefs}

\begin{nobs}
In the proof of the next lemma we use basic definitions from algebraic geometry
that can be found in \citet{MR0463157}.
For example the definitions of irreducible varieties, projective varieties, complete varieties, normal varieties,
invertible sheaves, Cartier divisors, Weyl divisors and basepoint free divisors.

Also, we use some definitions and concepts from toric geometry \citep{MR2810322}.
For example the toric variety associated to a fan, a torus-invariant divisor, a nef divisor
and finally, Demazure Vanishing.
We give a precise reference where the reader can find the definitions and results about toric geometry.
\end{nobs}

\begin{lema}\label{prel-mink}
Given a full dimensional lattice polytope $\mathcal{E}$, there exists an irreducible
projective normal toric variety $X$ such that every
$\N$-Minkowski summand $\mathcal{B}$ of $\mathcal{E}$ defines an invertible sheaf $\mathcal{O}_X(D)$
with
$$H^0(X,\mathcal{O}_X(D))=S_{\mathcal{B}},\quad H^p(X,\mathcal{O}_X(D))=0,\, p>0.$$
Even more, if $\mathcal{B}_1$ and $\mathcal{B}_2$ are two $\N$-Minkowski summands of $\mathcal{E}$
and $\mathcal{O}_X(D_1)$ and $\mathcal{O}_X(D_2)$ are the corresponding invertible sheaves
of $\mathcal{B}_1$ and $\mathcal{B}_2$ respectively,
then the invertible sheaf associated to $\mathcal{B}_1+\mathcal{B}_2$ is $\mathcal{O}_X(D_1+D_2)$.
\end{lema}
\begin{dem}
Given a full dimensional lattice polytope, we can construct a \emph{normal fan} $\Sigma$ \citep[Theorem 2.3.2]{MR2810322},
and a normal toric variety $X_{\Sigma}$ \citep[Theorem 3.1.5]{MR2810322}.

The normal fan associated to a full dimensional lattice polytope is \emph{complete} \citep[Proposition 2.3.8]{MR2810322}.
Then $X_{\Sigma}$ is also a complete variety \citep[Theorem 3.4.6]{MR2810322}.

There exists a more direct construction of $X_{\Sigma}$ using a multiple of the full dimensional
lattice polytope $\mathcal{E}$, but by Proposition 3.1.6 \citep{MR2810322} both constructions agree, $X_{\Sigma}\cong X_{\mathcal{E}}$.
The benefit of this direct construction is that $X_{\mathcal{E}}$ proves to be an irreducible projective variety.
Let us call $X$ the irreducible projective normal toric variety $X_{\Sigma}$.

Let $\mathcal{B}$ be an $\N$-Minkowski summand of $\mathcal{E}$. By Corollary 6.2.15 \citep{MR2810322}
there exists a torus invariant basepoint free Cartier divisor $D$ on $X$ such that
$$H^0(X,\mathcal{O}_X(D))=S_{\mathcal{B}}.$$
This last equality follows from Proposition 4.3.3 \citep{MR2810322} and the fact that
in a normal variety, every Cartier divisor is a Weyl divisor \citep[Definition 4.0.12]{MR2810322}.

Let us apply Demazure Vanishing, \citep[Theorem 9.2.3]{MR2810322}.
By definition, the \emph{support} of a complete fan is $\R^n$, \citep[Definition 3.1.18]{MR2810322}.
In particular, it has a \emph{convex support of full dimension}, \citep[\S 6.1]{MR2810322}.
Then the basepoint free Cartier divisor $D$ is \emph{nef}, \citep[Theorem 6.3.12]{MR2810322}.
Applying Demazure Vanishing, we obtain,
$$H^p(X,\mathcal{O}_X(D))=0,\quad p>0.$$

Let us prove the last paragraph of the lemma.
Let $D$ be a torus-invariant Cartier divisor on $X$.
Then there exists a polytope $\mathcal{P}_D$ such that
$H^0(X,\mathcal{O}_X(D))=S_{\mathcal{P}_D}$, \citep[Lemma, p. 66]{MR1234037}.
Even more, if $D$ is the torus-invariant basepoint free Cartier divisor associated to an $\N$-Minkowski
summand $\mathcal{B}$, then $\mathcal{P}_D=\mathcal{B}$, \citep[p. 68]{MR1234037}; \citep[Corollary 6.2.15]{MR2810322}.

Let $\mathcal{B}_1$ and $\mathcal{B}_2$ be two $\N$-Minkowski summands of $\mathcal{E}$
and $\mathcal{O}_X(D_1)$ and $\mathcal{O}_X(D_2)$ be the corresponding invertible sheaves
associated to $\mathcal{B}_1$ and $\mathcal{B}_2$ respectively.
Given that the sheaves
are generated by global sections, we
obtain $\mathcal{P}_{D_1+D_2}=\mathcal{P}_{D_1}+\mathcal{P}_{D_2}$ \citep[Exercise, p. 69]{MR1234037}.

Let $\mathcal{O}_X(D)$ be the invertible sheaf associated to the $\N$-Minkowski summand $\mathcal{B}_1+\mathcal{B}_2$
of $\mathcal{E}$. Then,
$$\mathcal{P}_D=\mathcal{B}_1+\mathcal{B}_2=\mathcal{P}_{D_1}+\mathcal{P}_{D_2}=\mathcal{P}_{D_1+D_2}.$$
This implies that $\mathcal{O}_X(D)\cong \mathcal{O}_X(D_1+D_2)$.
\hspace{\fill}\qed
\end{dem}

\begin{ndefs}
Let $\{f_1,\ldots,f_k\}$ be a sparse system in $(\C\setminus 0)^n$ with $r'<\infty$ solutions counted with multiplicities.
The torus $(\C\setminus 0)^n$ is contained in the variety $X$ of \autoref{prel-mink}
as an open subset, \citep[Definition 3.1.1]{MR2810322}. Homogenizing every equation
of the sparse system, we can consider the system in $X$.
For the homogenization process, see \citet[\S 5.4]{MR2810322}.

Let $Z\subseteq X$ be the zero-scheme of the resulting system and let $r\geq r'$ be the number of points in $Z$
counted with multiplicities.
We say that the sparse system $\{f_1,\ldots,f_k\}$
\emph{has no solution at infinity} if $r=r'$.
Otherwise, we say that it \emph{has solution at infinity}.
The number $r$ is called the \emph{expected number of solutions} of the sparse system.
\end{ndefs}

\section{Solving a sparse system.}\label{sec-spar}

The following notations and assumptions will be used in the rest of the section.
\begin{ass}\label{spar-ass}
Let $f_0,\ldots,f_k$ be Laurent polynomials
with Newton polytopes $\mathcal{A}_0,\ldots,\mathcal{A}_k$ respectively.
Let $\mathcal{B}_i=\mathcal{A}_0+\ldots+\widehat{\mathcal{A}_i}+\ldots+\mathcal{A}_k$, $0\leq i\leq k$
and let $\mathcal{E}=\mathcal{A}_0+\ldots+\mathcal{A}_k$.
Let $\mathcal{I}=\langle f_1,\ldots,f_k\rangle\subseteq S$ be the ideal generated by the sparse system.

Assume,
\begin{itemize}
\item $0\in \mathcal{A}_0$ and $f_0$ is a non-constant Laurent polynomial.
\item The lattice polytope $\mathcal{E}$ is full dimensional.
\item $\langle f_0,f_1,\ldots,f_k\rangle=S$.
\end{itemize}
\hspace{\fill}\qed
\end{ass}

\begin{nobs}
If $0\not\in\mathcal{A}_{0}$, we can divide the equation $f_{0}$ by some monomial
or we can consider the convex hull of $0$ and $\mathcal{A}_{0}$ as the new lattice polytope $\mathcal{A}_{0}$.
These
operations does not change the number of expected solutions nor the solutions in $(\C\setminus 0)^n$ of the sparse system.
Then without loss of generality, we can assume $0\in\mathcal{A}_0$. This assumptions implies
that $\mathcal{B}_0$ is contained in $\mathcal{E}=\mathcal{A}_0+\mathcal{B}_0$.

If $\mathcal{E}$ is not full dimensional, there exists an affine change of variables such that the
variables, say $x_{s+1},\ldots,x_n$, are missing in the sparse system.
This implies that we could work in $S=\C[x_1^{\pm 1},\ldots,x_s^{\pm 1}]$ making
$\mathcal{E}$ a full dimensional lattice polytope. This change
of variables involves the computation of Smith Normal Forms \citep{MR1135749}.
Another remark, is that it is easy to prove that if $\mathcal{A}_0$ is
full dimensional, then $\mathcal{E}$ is full dimensional.
Hence, we can consider $\mathcal{A}_0$ as a full dimensional lattice polytope.

It follows from $\langle f_0,f_1,\ldots,f_k\rangle=S$ that the associated zero-scheme in $X$ is empty.
We prove in the next theorem that a sparse system satisfying the previous assumptions
will have a finite number of expected solutions (possible zero).
This assumption is the most important one.
\end{nobs}

\begin{teo}\label{spar-teo1}
Same notation as before.
Suppose $f_0,\ldots,f_k$ are Laurent polynomials as in \autoref{spar-ass}. Then,
\begin{enumerate}[(a)]
\item The co-rank of the following linear map is the expected number of solutions (possibly zero),
$$\Phi:S_{\mathcal{B}_1}\times\ldots\times S_{\mathcal{B}_k}\rightarrow S_{\mathcal{E}},\quad
\Phi(g_1,\ldots,g_k)=\sum_{i=1}^k f_i\cdot g_i.$$
In particular, if the system has no solution at infinity, it is equal to the number of solutions in $(\C\setminus 0)^n$.
\item The lattice polytope $\mathcal{B}_0\subseteq \mathcal{E}$ satisfies
$$S_{\mathcal{B}_0}/
(\im(\Phi)\cap S_{\mathcal{B}_0})\cong S_{\mathcal{E}}/\im(\Phi).$$
\item The following linear map is surjective,
$$\Psi:S_{\mathcal{B}_0}\times S_{\mathcal{B}_1}\times\ldots\times S_{\mathcal{B}_k}\fl{} S_{\mathcal{E}},\quad
\Psi(g_0,g_1,\ldots,g_k)=f_0\cdot g_0+\sum_{i=1}^k f_i\cdot g_i.$$
\end{enumerate}
\end{teo}
\begin{dem}
Let us work with the projective variety $X$ of \autoref{prel-mink}.
For every integers $d_0,\ldots,d_k\geq 0$ consider the invertible sheaf $\mathcal{O}_X(d_0,\ldots,d_k)$
associated to the $\N$-Minkowski summand $d_0\mathcal{A}_0+\ldots+d_k\mathcal{A}_k$ of $\mathcal{E}$.
Then
$$H^0(X,\mathcal{O}_X(d_0,\ldots,d_k))=S_{d_0\mathcal{A}_0+\ldots+d_k\mathcal{A}_k},\quad
H^p(X,\mathcal{O}_X(d_0,\ldots,d_k))=0,\quad p>0.$$
Also, from the last paragraph of \autoref{prel-mink} we have the following property.
Let $d_i,d_i'$ be non-negative integers such that $d_i\geq d_i'\geq0$ for all $0\leq i\leq k$.
Then,
$$\mathcal{O}_X(d_0',\ldots,d_k')\otimes_{\mathcal{O}_X}\mathcal{O}_X(d_0-d_0',\ldots,d_k-d_k')
\cong\mathcal{O}_X(d_0,\ldots,d_k)\then$$
$$\mathcal{O}_X(d_0-d_0',\ldots,d_k-d_k')
\cong\mathcal{O}_X(d_0,\ldots,d_k)\otimes_{\mathcal{O}_X}\mathcal{O}_X(-d_0',\ldots,-d_k'),$$
where $\mathcal{O}_X(-d_0',\ldots,-d_k')$ denotes the dual sheaf of $\mathcal{O}_X(d_0',\ldots,d_k')$.

Let $e_i\in\Z^{k+1}$ be the vector with $1$ in the $(i+1)$-coordinate and $0$ in the rest, $0\leq i\leq k$.
For example, $e_0=(1,0,\ldots,0)$ and $e_{k}=(0,\ldots,0,1)$.
The Laurent polynomials $\{f_1,\ldots,f_k\}$ determine a $\mathcal{O}_X$-linear map
$$\mathcal{O}_X\rightarrow \mathcal{F},\quad \mathcal{F}=\mathcal{O}_X(e_1)\oplus\ldots\oplus\mathcal{O}_X(e_k),$$
given by $g\mapsto(f_1g,\ldots,f_kg)$.
Then, we can construct the dual Koszul complex associated to $\mathcal{F}$,
$$0\rightarrow \bigwedge^k\mathcal{F}^\vee\rightarrow\ldots\rightarrow\bigwedge^i\mathcal{F}^\vee\rightarrow\ldots
\rightarrow\mathcal{F}^\vee\rightarrow\mathcal{O}_X,$$
where $\mathcal{F}^\vee$ denotes the dual of $\mathcal{F}$ and
$$\bigwedge^s\mathcal{F}^\vee=\bigoplus_{1\leq i_1<\ldots<i_s\leq k}\mathcal{O}_X(-e_{i_1}-\ldots-e_{i_s}),\quad 2\leq s\leq k.$$

Let $Z\subseteq X$ be the zero scheme of the global section
$(f_1,\ldots,f_k)\in S_{\mathcal{A}_1}\oplus\ldots\oplus S_{\mathcal{A}_k}\cong H^0(X,\mathcal{F})$.
Let us prove that $Z$ is empty or of dimension $0$. Assume that $Z$ is not empty.
Let $H\subseteq X$ be the hypersurface given by the zeros of the
section $f_0\in S_{\mathcal{A}_0}\cong H^0(X,\mathcal{O}_X(1,0,\ldots,0))$.
Take an embedding of $X$ is some $\PP^N$ and let $\widehat{H}\subseteq\PP^N$ be an hypersurface
such that $\widehat{H}\cap X=H$. Given that the zero locus of $\{f_0,\ldots,f_k\}$ is empty in $X$,
we have $\emptyset=Z\cap H=Z\cap (\widehat{H}\cap X)=Z\cap \widehat{H}$.
Using Theorem 7.2 in \citet{MR0463157}, we obtain $\dim(Z)=0$.

Let us work with \emph{the augmented dual Koszul complex associated to $\mathcal{F}$},
$$0\rightarrow \bigwedge^k\mathcal{F}^\vee\rightarrow\ldots\rightarrow\bigwedge^i\mathcal{F}^\vee\rightarrow\ldots
\rightarrow\mathcal{F}^\vee\rightarrow\mathcal{O}_X\rightarrow\mathcal{O}_Z\rightarrow 0.$$
By \S 2, 1B, Proposition 1.4 (a) \citep{MR1264417} it is an exact complex.

Let $U\subseteq X$ be an affine
open subset containing $Z$. Let $T$ be the coordinate ring of $U$ and let $\mathcal{J}$ be the
ideal of $Z\subseteq U$. Then,
$$H^0(X,\mathcal{O}_Z(d_0,\ldots,d_k))=H^0(U,\mathcal{O}_Z)=T/\mathcal{J},\quad \forall d_0,\ldots,d_k\geq 0.$$

Recall from Proposition 2.9 \citep{MR0463157} that cohomology commutes with direct sums and from
Theorem 6.0.18 and Proposition 6.0.17 \citep{MR2810322} that invertible sheaves are locally free.

\begin{enumerate}[(a)]
\item The exactness of the augmented dual Koszul complex associated to $\mathcal{F}$
is preserved by twisting with the invertible sheaf $\mathcal{O}_X(1,\ldots,1)$,
and giving that each term of the resulting complex has no higher cohomology,
the following complex of vector spaces is exact \citep[\S 2, 2A, Lemma 2.4]{MR1264417},
$$S_{\mathcal{B}_1}\times\ldots\times S_{\mathcal{B}_k}\fl{\Phi}S_{\mathcal{E}}\rightarrow T/\mathcal{J}\rightarrow 0.$$
If the sparse system has no solution at infinity, we can take the torus as the open set $U$,
then $Z\subseteq (\C\setminus 0)^{n}$ and $T/\mathcal{J}=S/\mathcal{I}$.
\item In a similar way, twisting the augmented dual Koszul complex associated to $\mathcal{F}$ with the
invertible sheaf $\mathcal{O}_X(0,1,\ldots,1)$,
the following map is surjective,
$$S_{\mathcal{B}_0}\rightarrow T/\mathcal{J}\rightarrow 0.$$
Let $K$ be the kernel of $S_{\mathcal{B}_0}\rightarrow T/\mathcal{J}$. Then,
$$\xymatrix{
0\ar[r]&K\ar@{^{(}->}[r]\ar@{^{(}-->}[d]&S_{\mathcal{B}_0}\ar@{}[dr]|{\equiv}\ar@{^{(}->}[d]\ar[r]&T/\mathcal{J}\ar@{=}[d]\ar[r]&0\\
0\ar[r]&\im(\Phi)\ar@{^{(}->}[r]&S_{\mathcal{E}}\ar[r]&T/\mathcal{J}\ar[r]&0}$$
The inclusion $\mathcal{B}_0\subseteq\mathcal{E}$ follows from the assumption $0\in\mathcal{A}_0$.
Given that both rows are exact, $K$ must be equal to $S_{\mathcal{B}_0}\cap\im(\Phi)$. Then,
$$S_{\mathcal{B}_0}/(S_{\mathcal{B}_0}\cap\im(\Phi))\cong T/\mathcal{J}\cong S_{\mathcal{E}}/\im(\Phi).$$
\item  Finally, given that the zero locus of $\{f_0,\ldots,f_k\}$ is empty in $X$, we can use similar arguments
as before with the sheaf $\mathcal{F}'=\mathcal{O}_X(e_0)\oplus\ldots\oplus\mathcal{O}_X(e_k)$
to prove that the following map is surjective,
$$S_{\mathcal{B}_0}\times\ldots\times S_{\mathcal{B}_k}\fl{\Psi}S_{\mathcal{E}}\rightarrow 0.$$
Specifically, the augmented dual Koszul complex associated to $\mathcal{F}'$ is
$$0\rightarrow \bigwedge^{k+1}\mathcal{F}'^\vee\rightarrow\ldots\rightarrow\bigwedge^i\mathcal{F}'^\vee\rightarrow\ldots
\rightarrow\mathcal{F}'^\vee\rightarrow\mathcal{O}_X\rightarrow 0.$$
The result follows by twisting it with $\mathcal{O}_X(1,\ldots,1)$ and taking global sections.
\end{enumerate}
\hspace{\fill}\qed
\end{dem}

\begin{nobs}
From the previous proof, part (c), we obtain a formula involving the number of lattice points in $\mathcal{E}$.
The augmented dual Koszul complex associated to
$\mathcal{F}'=\mathcal{O}_X(e_0)\oplus\ldots\oplus\mathcal{O}_X(e_k)$
twisted by $\mathcal{O}_X(1,\ldots,1)$ is exact
and each term has no higher cohomology.
Hence its Euler characteristic is zero,
$$\#(\mathcal{E}\cap\Z^n)-\sum_{i=0}^k\#( (\mathcal{A}_0+\ldots+\widehat{\mathcal{A}_i}+\ldots+\mathcal{A}_k)\cap\Z^n)+\ldots
-(-1)^k \sum_{i=0}^k \#(\mathcal{A}_i\cap\Z^n)+(-1)^k=0.$$
This formula is similar to the alternate volume formula in \citet{MR0435072}.

For example, consider the simplex $\Delta\subseteq\R^3$, $\Delta=\mathcal{A}_0+\mathcal{A}_1+\mathcal{A}_2$,
where
$$\mathcal{A}_0=\text{conv}((0,0,0),(p,0,0)),\quad
\mathcal{A}_1=\text{conv}((0,0,0),(0,q,0)),\quad \mathcal{A}_2=\text{conv}((0,0,0),(0,0,r))$$
and $p,q,r$ are three positive prime numbers. Then $\#(\Delta\cap\Z^3)$ is equal to
$$(p+q+1)+(p+r+1)+(q+r+1)-(p+1)-(q+1)-(r+1)+1=p+q+r+1.$$
When $p=q=r=1$, the standard simplex in $\R^3$ has $4$ points in $\Z^3$.

For more on counting points in a lattice polytope, see \citet{MR2159956}.
\hspace{\fill}\qed
\end{nobs}

The following corollary is an adaptation of a theorem in \citet{MR1577000}.
\begin{cor}\label{spar-null}
Suppose $f_0,\ldots,f_k$ are Laurent polynomials
with Newton polytopes $\mathcal{A}_0,\ldots,\mathcal{A}_k$ respectively.
Let $\mathcal{B}_i=\mathcal{A}_0+\ldots+\widehat{\mathcal{A}_i}+\ldots+\mathcal{A}_k$, $0\leq i\leq k$
and let $\mathcal{E}=\mathcal{A}_0+\ldots+\mathcal{A}_k$.

Let $\Psi$ be the following linear map,
$$\Psi:S_{\mathcal{B}_0}\times S_{\mathcal{B}_1}\times\ldots\times S_{\mathcal{B}_k}\fl{} S_{\mathcal{E}},\quad
\Psi(g_0,g_1,\ldots,g_k)=f_0\cdot g_0+\sum_{i=1}^k f_i\cdot g_i.$$
Assume that $\mathcal{E}$ is full dimensional. Then,
$$\rk(\Psi)=\#(\mathcal{E}\cap\Z^n)\iff \langle f_0,\ldots,f_k\rangle=S.$$
\end{cor}
\begin{dem}
If $\langle f_0,\ldots,f_k\rangle=S$, by \autoref{spar-teo1}(c), $\Psi$ is surjective.
Analogously, if $\Psi$ is surjective, there exists a monomial
$x^m\in S_{\mathcal{E}}=\im(\Psi)\subseteq \langle f_0,\ldots,f_k\rangle$. In particular, $1\in\langle f_0,\ldots,f_k\rangle$.
\hspace{\fill}\qed
\end{dem}

\begin{nnotc}
Notations and assumptions as in \autoref{spar-ass}.
Consider $\mathcal{E}\cap\Z^n$ and $\mathcal{B}_i\cap\Z^n$ as
monomial ordered bases of $S_{\mathcal{E}}$ and $S_{\mathcal{B}_i}$ respectively, $0\leq i\leq k$.
Let $p$ be the cardinal of $\mathcal{B}_0\cap\Z^n$, $p_i$ the cardinal of $\mathcal{B}_i\cap\Z^n$, $1\leq i\leq k$
and $p+q$ be the cardinal of $\mathcal{E}\cap\Z^n$.
Given that $f_0$ is not a constant in $S_{\mathcal{A}_0}$,
the inclusion $\mathcal{B}_0\subseteq\mathcal{E}$ is proper, thus $q>0$.
$$\mathcal{B}_0\cap\Z^n=\{m_1,\ldots,m_p\},\quad
\mathcal{E}\cap\Z^n=\{m_1,\ldots,m_p,m_{p+1},\ldots,m_{p+q}\},$$
where $m_i$ is a point in $\Z^n$, $1\leq i\leq p+q$.

Let us define the coefficient matrix $M$ associated to the sparse system $\{f_1,\ldots,f_k\}$ and $f_0$.
Let $M\in\C^{(p+q)\times (p+p_1+\ldots+p_k)}$ be the rectangular matrix associated to $\Psi$ in these bases,
$$M=\begin{pmatrix}
M_{11}&M_{12}\\
M_{21}&M_{22}
\end{pmatrix},\quad M_{11}\in\C^{p\times p},\,M_{22}\in\C^{q\times (p_1+\ldots+p_k)}.$$
Then,
$$\begin{pmatrix}
x^{m_1}\ldots x^{m_p}&x^{m_{p+1}}\ldots x^{m_{p+q}}
\end{pmatrix}
\begin{pmatrix}
M_{11}&M_{12}\\
M_{21}&M_{22}
\end{pmatrix}=
\begin{pmatrix}
f_0x^{m_1}\ldots f_0x^{m_p}&f_1\cdot\mathcal{B}_1\ldots f_k\cdot\mathcal{B}_k
\end{pmatrix},$$
where $f_i\cdot\mathcal{B}_i$ is the row vector obtained by multiplying $f_i$ with the monomials
in $\mathcal{B}_i\cap\Z^n$, $1\leq i\leq k$.
We are abusing the notation;
the point $m=(m_1,\ldots,m_n)\in\Z^n$ corresponds to the monomial $x^m=x_1^{m_1}\ldots x_n^{m_n}$.
\hspace{\fill}\qed
\end{nnotc}

\begin{teo}\label{spar-F}
Same notation as before.
Suppose $f_0,\ldots,f_k$ are Laurent polynomials as in \autoref{spar-ass}.
Let $F\in\C^{(p_1+\ldots+p_k)\times p}$ be a solution
of the linear equation $M_{21}+M_{22}F=0$.

Then, every solution, $\xi$, of the sparse system determines a
left eigenvector of $M_{11}+M_{12}F$. Even more,
$f_0(\xi)$ is the eigenvalue of that left eigenvector.
\end{teo}
\begin{dem}
Let us see that the hypotheses given in \autoref{spar-ass} imply that the rectangular matrix $M_{22}$ has full rank $q$.
The matrix $M_{22}$ is the matrix of the composition of the following maps,
$$S_{\mathcal{B}_1}\times\ldots\times S_{\mathcal{B}_k}\fl{\Phi} S_{\mathcal{E}}\fl{\pi}S_{\mathcal{E}}/S_{\mathcal{B}_0},$$
where $\pi$ is the quotient map. Then the rank of $M_{22}$, $t$, is equal to $\rk(\pi\,\Phi)$ and
$$\im(\pi\Phi)=\pi(\im(\Phi))=
\im(\Phi)/\left(\im(\Phi)\cap S_{\mathcal{B}_0}\right)\then$$
$$t:=\rk(M_{22})=\dim(\im(\Phi)/\left(\im(\Phi)\cap S_{\mathcal{B}_0}\right))=\dim(\im(\Phi))-\dim(\im(\Phi)\cap S_{\mathcal{B}_0}).$$
Note that the dimension of $S_{\mathcal{E}}/S_{\mathcal{B}_0}$ is equal to $q$,
$$\dim(S_{\mathcal{E}}/S_{\mathcal{B}_0})=\dim(S_{\mathcal{E}})-\dim(S_{\mathcal{B}_0})=(p+q)-p=q.$$
Let us prove $q=t$. Using \autoref{spar-teo1}(b), we get
$$\dim(S_{\mathcal{E}}/\im(\Phi))=\dim(S_{\mathcal{B}_0}/\left(\im(\Phi)\cap S_{\mathcal{B}_0}\right))\then $$
$$\dim(S_{\mathcal{E}})-\dim(\im(\Phi))=\dim(S_{\mathcal{B}_0})-\dim(\im(\Phi)\cap S_{\mathcal{B}_0})\then$$
$$q=\dim(S_{\mathcal{E}})-\dim(S_{\mathcal{B}_0})=\dim(\im(\Phi))-\dim(\im(\Phi)\cap S_{\mathcal{B}_0})=t.$$

Now that we know that $\rk(M_{22})=q$, it is easy to prove that
there exists a matrix
$F$ such that
$$M_{22}F=-M_{21},\quad F\in\C^{(p_1+\ldots+p_k)\times p}.$$
Each column of $F$, $c_1,\ldots,c_p$, is a solution of the linear system $M_{22}c_i=b_i$,
where $b_i\in\C^q$ is the $i$-column vector of $-M_{21}$, $1\leq i\leq p$.

Let $\xi\in\C^n$ be a solution of the sparse system, $f_1(\xi)=\ldots=f_k(\xi)=0$. Then
{\small
$$\begin{pmatrix}
\xi^{m_1}\ldots \xi^{m_p}&\xi^{m_{p+1}}\ldots \xi^{m_{p+q}}
\end{pmatrix}
\begin{pmatrix}
M_{11}&M_{12}\\
M_{21}&M_{22}
\end{pmatrix}=
f_0(\xi)\cdot\begin{pmatrix}
\xi^{m_1}\ldots \xi^{m_p}&0
\end{pmatrix}\then$$
$$\begin{pmatrix}
\xi^{m_1}\ldots \xi^{m_p}&\xi^{m_{p+1}}\ldots \xi^{m_{p+q}}
\end{pmatrix}
\begin{pmatrix}
M_{11}&M_{12}\\
M_{21}&M_{22}
\end{pmatrix}
\begin{pmatrix}
I&0\\
F&I
\end{pmatrix}=
f_0(\xi)\cdot\begin{pmatrix}
\xi^{m_1}\ldots \xi^{m_p}&0
\end{pmatrix}
\begin{pmatrix}
I&0\\
F&I
\end{pmatrix}\then$$
$$(\xi^{m_1}\ldots \xi^{m_p})\left(M_{11}+M_{12}F\right)=f_0(\xi)\cdot (\xi^{m_1} \ldots \xi^{m_p}).$$
}Then $\xi$ determines a left eigenvector of $M_{11}+M_{12}F$ with eigenvalue $f_0(\xi)$.
\hspace{\fill}\qed
\end{dem}

\begin{nobs}
In the previous theorem we proved that every solution of a sparse system as in \autoref{spar-ass}
determines a left eigenvector of the square matrix $M_{11}+M_{12}F$. If the
geometric multiplicity of an eigenvalue (the dimension of its left eigenspace) is greater than one,
then we cannot use the computation of left eigenvectors to deduce the solutions of
the sparse system.
\hspace{\fill}\qed
\end{nobs}

Let us relate the multiplicity of a root $\xi$ with the multiplicity of the eigenvalue
$f_0(\xi)$ in the matrix $M_{11}+M_{12}F$.
\begin{prop}\label{spar-mult}
Same notation as before.
Suppose $f_0,\ldots,f_k$ are Laurent polynomials as in \autoref{spar-ass}.
Let $\xi\in(\C\setminus 0)^n$ be a solution of the sparse system with multiplicity $\mu$.
Then, the eigenvalue $f_0(\xi)$ of $M_{11}+M_{12}F$ has multiplicity
greater than or equal to $\mu$.
\end{prop}
\begin{dem}
The characteristic polynomial
of the multiplication map $M_{f_0}:S/\mathcal{I}\rightarrow S/\mathcal{I}$ is
$$\chi(t)=(t-f_0(\xi_1))^{\mu_1}\ldots(t-f_0(\xi_s))^{\mu_s},$$
where $\xi_1,\ldots,\xi_s$ are the solutions of the sparse system in $(\C\setminus 0)^n$ and $\mu_1,\ldots,\mu_s$
their respective multiplicities (see \citet[2.1.14]{MR2161984}).

Let us relate the multiplication matrix $M_{f_0}$ with our matrix $M$.
Recall that
the columns of the matrix of $\Phi$ are the multiples of $\{f_1,\ldots,f_k\}$,
$$[\Phi]=\begin{pmatrix}
M_{12}\\
M_{22}
\end{pmatrix}.$$
Also, that the matrix of the map $\Psi|_{S_{\mathcal{B}_0}}:S_{\mathcal{B}_0}\rightarrow S_{\mathcal{E}}$
corresponds to the multiples of $f_0$,
$$\begin{pmatrix}
M_{11}\\
M_{21}
\end{pmatrix}.$$
In order to find the class of $f_0\in S/\mathcal{I}$, we need to add/substract monomial multiples of $\{f_1,\ldots,f_k\}$
to $f_0$. This process may be done by column operations in $M$. In particular,
the matrix
$$\begin{pmatrix}
M_{11}+M_{12}F&M_{12}\\
0&M_{22}
\end{pmatrix}=
\begin{pmatrix}
M_{11}&M_{12}\\
M_{21}&M_{22}
\end{pmatrix}
\begin{pmatrix}
I&0\\
F&I
\end{pmatrix}$$
obtained by column operations from $M$, also determines the class in $S/\mathcal{I}$ of $f_0$.
Specifically, the class of $x^{m_j} f_0$ in $S/\mathcal{I}$
is the same as the class of the $j$-column of $M$ and the $j$-column of $M_{11}+M_{12}F$, $1\leq j\leq p$,
$$x^{m_j} f_0\equiv\sum_{i=1}^{p+q}x^{m_i}a_{ij} \equiv\sum_{i=1}^{p}x^{m_i}b_{ij}\quad\text{mod }\mathcal{I},$$
where $a_{ij}=M_{ij}$, $b_{ij}=(M_{11}+M_{12}F)_{ij}$ and
the first $p$ monomials are in $\mathcal{B}_0\subseteq \mathcal{E}$ and the last $q$ monomials are in $\mathcal{E}\setminus\mathcal{B}_0$.

Let us call $\sigma_{f_0}$ the map associated to $M_{11}+M_{12}F$,
$$\sigma_{f_0}:S_{\mathcal{B}_0}\rightarrow S_{\mathcal{B}_0},\quad
\sigma_{f_0}(x^{m_j})=\sum_{i=1}^p x^{m_i}b_{ij}\quad 1\leq j\leq p.$$
Then,
we have the following commutative diagram,
$$\xymatrix{
S_{\mathcal{B}_0}\ar@{}[dr]|{\equiv}\ar[d]_{\sigma_{f_0}}\ar[r]^\pi&S/\mathcal{I}\ar[d]^{M_{f_0}}\\
S_{\mathcal{B}_0}\ar[r]^\pi&S/\mathcal{I}}$$

By \autoref{spar-teo1}(b), the map $\pi$ above is an epimorphism.
Hence the characteristic polynomial of $M_{f_0}$ divides the characteristic polynomial of the matrix $M_{11}+M_{12}F$,
that is, $\chi_{\sigma_{f_0}}(t)=\chi(t)P(t)$, where $P(t)$ is some polynomial (it depends on $F$).
\hspace{\fill}\qed
\end{dem}

To end this section, let us give an example on how to apply the previous theorems,
\begin{nejemplo}
Consider the intersection of a line with a parabola,
$$\left\{\begin{array}{rcccl}
f_1&=&1+x+y&=&0\\
f_2&=&1+x^2+y&=&0
\end{array}\right.$$
They intersect in $(1,-2)$ and $(0,-1)\in\C^2$.
Note that the ideal generated by $\langle f_1,f_2\rangle$ is radical.
Let $f_0=x-2y$ be a linear form. The value of $f_0$ at each solution is $f_0(1,-2)=5$ and $f_0(0,-1)=2$.

Let us identify the monomial $x^ny^m$ with the point $(n,m)\in\Z^2$.
Take the lattice polytopes associated to $f_0,f_1$ and $f_2$,
$$\mathcal{A}_0\cap\Z^2=\mathcal{A}_1\cap\Z^2=\{1,x,y\}, \quad
\mathcal{A}_2\cap\Z^2=\{1,y,x,x^2\}.$$
Then,
$$\mathcal{B}_0\cap\Z^2=\{1, x^2y, x^2, y^2, x^3, xy, x, y\},\quad
\mathcal{E}\cap\Z^2=\{1, x^2y, x^2, y^2, x^3, xy, x, y, x^3y, xy^2, y^2x^2, x^4, y^3\},$$
In these bases the coefficient matrix $M$ is equal to,
{\footnotesize
$$\left(\begin{array}{cccccccc|cccccccccccccc}
0&0&0&0&0&0&0&0&1&0&0&0&0&0&0&0&1&0&0&0&0&0\\
0&0&1&0&-2&0&0&0&0&1&1&0&1&0&0&0&0&0&0&1&0&1\\
0&0&0&0&0&0&-2&1&0&0&1&0&0&0&1&1&0&1&0&0&1&0\\
0&0&0&0&0&0&0&-2&0&0&0&1&0&0&0&1&0&0&1&0&0&1\\
0&0&0&0&0&0&1&0&0&0&0&0&1&0&1&0&1&0&0&1&0&0\\
0&0&0&0&1&0&0&0&0&0&0&0&1&1&0&0&0&0&0&0&1&0\\
1&0&0&0&0&0&0&0&1&0&0&0&0&0&1&0&0&0&0&0&1&0\\
-2&0&0&0&0&0&0&0&1&0&0&0&0&0&0&1&1&0&0&0&0&1\\
\hline
0&1&0&0&0&-2&0&0&0&1&0&0&0&1&0&0&0&1&0&0&0&0\\
0&0&-2&1&0&0&0&0&0&0&1&1&0&0&0&0&0&1&0&0&0&0\\
0&-2&0&0&0&0&0&0&0&1&0&0&0&0&0&0&0&0&1&0&0&0\\
0&0&0&0&0&1&0&0&0&0&0&0&0&1&0&0&0&0&0&1&0&0\\
0&0&0&-2&0&0&0&0&0&0&0&1&0&0&0&0&0&0&1&0&0&0
\end{array}\right).$$
}
Using \autoref{spar-null},
$\rk(M)=13=\#(\mathcal{E}\cap\Z^2)$ implies that
the ideal $\langle f_0,f_1,f_2\rangle$ is $S$ as we already knew.
Also, we can recover the expected number of solutions, \autoref{spar-teo1}(a),
$$\#(\mathcal{E}\cap\Z^2)-
\rk\begin{pmatrix}
M_{12}\\
M_{22}
\end{pmatrix}=13-11=2.$$
Let us compute a matrix $F$ and $M_{11}+M_{12}F$,
{\footnotesize
$$F=\left(\begin{array}{cccccccc}
0&0&0&0&0&0&0&0\\0&2&0&0&0&0&0&0\\0&3&2&-3&0&-3&0&0\\
0&0&0&2&0&0&0&0\\0&0&0&0&0&0&0&0\\0&0&0&0&0&-1&0&0\\
0&0&0&0&0&0&0&0\\0&0&0&0&0&0&0&0\\0&0&0&0&0&0&0&0\\
0&-3&0&0&0&3&0&0\\0&0&0&0&0&0&0&0\\0&0&0&0&0&0&0&0\\
0&0&0&0&0&0&0&0\\0&0&0&0&0&0&0&0
\end{array}\right),\quad
M_{11}+M_{12}F=
\left(\begin{array}{cccccccc}
0&0&0&0&0&0&0&0\\
0&5&3&-3&-2&-3&0&0\\
0&0&2&-3&0&0&-2&1\\
0&0&0&2&0&0&0&-2\\
0&0&0&0&0&0&1&0\\
0&0&0&0&1&-1&0&0\\
1&0&0&0&0&0&0&0\\
-2&0&0&0&0&0&0&0
\end{array}\right).$$
}

The characteristic polynomial of $M_{11}+M_{12}F$ is equal to $t^4(t+1)(t-2)^2(t-5)$
and its minimal polynomial is $t^3(t+1)(t-2)^2(t-5)$.
The left eigenspace associated to $2$ is $\langle (1, 0, 0, 1, 0, 0, 0, -1)\rangle$,
and the left eigenspace associated to $5$ is $\langle (1, -2, -2, 4, 1, 1, 1, -2)\rangle$.
Then, looking at the last two coordinates (the monomials $x$ and $y$),
we get the two solutions $(0,-1)$ and $(1,-2)$.

Note that the characteristic polynomial and the minimal polynomial of $M_{11}+M_{12}F$ are different.
Also, that the multiplicity of $2=f_0(0,-1)$ is two.

In this example, it is possible to choose another $F$ (without changing $f_0$)
to get eigenvalues with multiplicity one in $M_{11}+M_{12}F$,
{\footnotesize
$$F'=\left(\begin{array}{cccccccc}
0&3&-2&-9&0&1&0&0\\
0&1&4&8&0&-6&0&0\\
0&1&0&-4&0&3&0&0\\
0&-1&4&10&0&-6&0&0\\
0&7&-8&-28&0&15&0&0\\
0&-2&-2&-1&0&5&0&0\\
0&-3&-2&1&0&2&0&0\\
0&-5&4&16&0&-5&0&0\\
0&-1&-2&-3&0&1&0&0\\
0&0&-2&-7&0&3&0&0\\
0&1&-4&-8&0&6&0&0\\
0&2&2&1&0&-6&0&0\\
0&-5&6&21&0&-10&0&0\\
0&2&0&-2&0&1&0&0
\end{array}\right).$$
}
\hspace{\fill}\qed
\end{nejemplo}

\begin{nobs}[About the size of the matrices]
It is important to mention that our matrix construction produce an extremely large matrix.
In the previous example we produced a matrix in $\C^{13\times 22}$,
but considering different monomial bases it is possible to construct a smaller matrix in $\C^{6\times 7}$.
The following monomials were suggested by a referee.
Let $\mathcal{B}_0\cap\Z^2=\mathcal{B}_1\cap\Z^2=\{1,x,y\}$ and $\mathcal{B}_2\cap\Z^2=\{1\}$. Then, the coefficient matrix
in these bases is
{\footnotesize
$$M=\left(\begin{array}{ccc|cccc}
0&0&0&1&0&0&1\\
1&0&0&1&1&0&0\\
-2&0&0&1&0&1&1\\
\hline
0&1&0&0&1&0&1\\
0&-2&1&0&1&1&0\\
0&0&-2&0&0&1&0
\end{array}\right)$$
}
Applying the same procedure as before, the matrix $M_{11}+M_{12}F$ has three eigenvalues, $0,2$ and $5$
with left eigenvectors $(1,6,4)$, $(1,0,-1)$ and $(1,1,-2)$ respectively.

Let us explain briefly why this matrix works. In this example, the polytope $\mathcal{E}=\mathcal{A}_0+\mathcal{B}_0$
satisfies $\mathcal{E}\cap\Z^2=\{1,x,y,x^2,xy,y^2\}$, then the projective variety $X$ of \autoref{prel-mink}
is the projective plane $\PP^2$. The invertible sheaf associated to $\mathcal{O}_X(d_0,d_1,d_2)$
is equal to $\mathcal{O}_{\PP^2}(d_0+d_1+2d_2)$.
Then, the augmented dual Koszul complex associated to the section $(f_1,f_2)$ is
$$0\rightarrow\mathcal{O}_{\PP^2}(-3)\rightarrow\mathcal{O}_{\PP^2}(-1)\oplus
\mathcal{O}_{\PP^2}(-2)\rightarrow\mathcal{O}_{\PP^2}\rightarrow\mathcal{O}_{Z}\rightarrow 0.$$
Using Theorem III.5.1 in \citet{MR0463157}, we know that the shaves $\mathcal{O}_{\PP^2}(-1)$ and $\mathcal{O}_{\PP^2}(-2)$
have no higher cohomology.

In order to prove that $M_{22}$ has full rank we need to prove \autoref{spar-teo1}(b)
and use the first part of the proof in \autoref{spar-F}.
Twisting the previous complex by $\mathcal{O}_{\PP^2}(1)$ and taking global sections,
we obtain that the map $S_{\mathcal{B}_0}\rightarrow S/\mathcal{I}$ is surjective.
Hence, $M_{22}$ has full rank.

To prove that the co-rank of the column block matrix $(M_{12};M_{22})$ is $2$,
we need to prove \autoref{spar-teo1}(a).
It follows by twisting with $\mathcal{O}_{\PP^2}(2)$ and taking global sections.

In the same way, it is easy to prove that $M$ has full rank
twisting by $\mathcal{O}_{\PP^2}(2)$ and
taking global sections the dual Koszul
complex associated to $(f_0,f_1,f_2)$
$$0\rightarrow\mathcal{O}_{\PP^2}(-4)\rightarrow
\mathcal{O}_{\PP^2}(-2)\oplus\mathcal{O}_{\PP^2}(-2)\oplus\mathcal{O}_{\PP^2}(-3)
\rightarrow\mathcal{O}_{\PP^2}(-1)\oplus\mathcal{O}_{\PP^2}(-1)\oplus
\mathcal{O}_{\PP^2}(-2)\rightarrow\mathcal{O}_{\PP^2}\rightarrow 0.$$

All of our results are based on Demazure Vanishing in \autoref{prel-mink}. In this
particular example the sheaves $\mathcal{O}_{\PP^2}(-1)$ and $\mathcal{O}_{\PP^2}(-2)$
have no higher cohomology by different reasons. It is worth mentioning that the construction of $M$ can be improved
to produce smaller matrices.
The size is controlled by the lattice polytope $\mathcal{E}$
and the space of global sections $H^0(X,\mathcal{F}'^\vee(1,\ldots,1))$,
where $X=X_{\mathcal{E}}$ is the projective variety of \autoref{prel-mink}
and $\mathcal{F}'$ is the sheaf associated to the sparse system and $f_0$.
Our approach gives a \emph{canonical} matrix to work in general.

Let us mention the following related result on reducing
the size of the coefficient matrix.
In \citet{MR1251972}, the authors work with (essentially) the same map $\Psi$, the first Koszul map,
to produce a formula of the sparse resultant of $n+1$ Laurent polynomials in $n$ variables.
Using a \emph{Row content function} they constructed a square coefficient matrix such that
its determinant is a non-zero multiple of the sparse resultant.
Continuing this work, in \citet{MR1374227}, the authors proposed an
incremental algorithm to obtain this submatrix.
Finally, in \citet{MR2004032}, the authors
provided a coefficient matrix of optimal size for the case of multihomogeneous systems.
\end{nobs}

\section{Application.}\label{sec-app}
To conclude this article, let us give an application to approximate the
maximum of a generic trilinear form over a product of spheres,
$$\ell:\R^{n+1}\times\R^{m+1}\times\R^{s+1}\rightarrow\R,\quad \ell(x,y,z)=\sum_{(i,j,k)=0}^{(n,m,s)}a_{ijk}x_iy_jz_k,
\quad \max_{\|x\|=\|y\|=\|z\|=1}|\ell(x,y,z)|,$$
where the norm is the usual 2-norm.
This problem was studied in \citet{MR3003946}.
In the literature, the maximum of $\ell$ over a product of spheres,
is called the first singular value of $\ell$ \citep[\S 3]{1574201}.

Using Lagrange method of multipliers \citep[\S 13.7]{MR0344384} the extreme points of $\ell$
over a product of spheres, $\mathbb{S}^{n}\times\mathbb{S}^{m}\times\mathbb{S}^{s}$, satisfy
$$\left\{\begin{array}{lcr}
\partial\ell/\partial x_i(x_0,\ldots,x_n,y_0,\ldots,y_m,z_0,\ldots,z_s)&=&2\alpha x_i,\quad 0\leq i\leq n,\\
\partial\ell/\partial y_j(x_0,\ldots,x_n,y_0,\ldots,y_m,z_0,\ldots,z_s)&=&2\beta y_j,\quad 0\leq j\leq m,\\
\partial\ell/\partial z_k(x_0,\ldots,x_n,y_0,\ldots,y_m,z_0,\ldots,z_s)&=&2\lambda z_k,\quad 0\leq k\leq s,
\end{array}\right.$$
$$\quad\alpha,\beta,\lambda\in\R,\quad\|x\|=\|y\|=\|z\|=1.$$
These equations imply that the vector $\partial\ell/\partial x(x,y,z)$ is a multiple of $x$. Same for $y$ and $z$.
In other words, considering the system in $\PP^n\times \PP^m\times\PP^s$, we can hide the variables $\alpha,\beta$ and $\lambda$,
$$\left\{\begin{array}{lcr}
x_j\partial\ell/\partial x_i(x,y,z)&=&x_i\partial\ell/\partial x_j(x,y,z),\quad 0\leq i<j\leq n,\\
y_j\partial\ell/\partial y_i(x,y,z)&=&y_i\partial\ell/\partial y_j(x,y,z),\quad 0\leq i<j\leq m,\\
z_j\partial\ell/\partial z_i(x,y,z)&=&z_i\partial\ell/\partial z_j(x,y,z),\quad 0\leq i<j\leq s.
\end{array}\right.$$
Given that $\ell$ is trilinear, the expression $x_j\partial\ell/\partial x_i(x,y,z)$ is equal
to $\ell(x_je_i,y,z)$, where $e_i\in\R^{n+1}$ is the vector with $1$ in the $i$-coordinate and $0$ in the rest.
Same for $y$ and $z$.
Summing up, the extreme points of $\ell$ satisfy
the following system of equations in $\R^{n+1}\times\R^{m+1}\times\R^{s+1}$,
$$\left\{\begin{array}{rcl}
\ell(x_je_i-x_ie_j,y,z)&=&0,\quad 0\leq i<j\leq n,\\
\ell(x,y_je_i-y_ie_j,z)&=&0,\quad 0\leq i<j\leq m,\\
\ell(x,y,z_je_i-z_ie_j)&=&0,\quad 0\leq i<j\leq s,\\
x_0^2+\ldots+x_n^2&=&1\\
y_0^2+\ldots+y_m^2&=&1\\
z_0^2+\ldots+z_s^2&=&1\\
\end{array}\right.$$

Assume that $\ell$ is generic and $2n,2m,2s\leq n+m+s$ \citep[\S 14, 1.3]{MR1264417}, hence
the extreme points are finite with multiplicity one.
Enumerate the equations, $f_1,\ldots,f_k$.
Let $\lambda_1,\ldots,\lambda_r$ be the real eigenvalues of $M_{11}+M_{12}F$ associated to solutions of
the system $f_1=\ldots=f_k=0$ and to the generic trilinear form $\ell$.
If $|\lambda_1|\geq|\lambda_i|$, $2\leq i\leq r$,
then $|\lambda_1|$ is the maximum value of $\ell$ over $\mathbb{S}^{n}\times\mathbb{S}^{m}\times\mathbb{S}^{s}$.

\begin{nobs}[Numerical Issues]
The previous application has several numerical issues.
For example, the computation of eigenvalues \citep{MR1810531,MR1810528,MR1906741,MR1809979,MR2043496},
the test of the genericity of $\ell$ and also, the evaluation of
a possible root in the equations. These last issues, may be solved using \emph{interval arithmetics}, \citep{MR1708669,MR2863516}.

In this application we did not required the computation of eigenvectors.
It is a delicate numerical issue.
When the matrix $M_{11}+M_{12}F$ is non-derogatory,
we can apply the work in \citet{MR2068273}. See also, \citet{MR1799632,MR1810528,MR684183,MR1318956,MR1213180}.
In other cases, it is possible to adapt some ideas from \citet{MR1330867}.
For a method that works on general systems, we refer to \citet{MR2742705}.
\hspace{\fill}\qed
\end{nobs}

\section*{Acknowledgments.}

This work was supported by CONICET, Argentina.
The author would like to thank
Alicia Dickenstein for reading the first version of this paper.
Thanks are also due to the anonymous referees for their useful suggestions and ideas.

\section*{References.}

\def\cprime{$'$} \def\cprime{$'$}

\end{document}